\documentclass[12pt]{article}

\usepackage{mathrsfs}
\usepackage{amsfonts}
\usepackage{amssymb,amsmath,amsthm, amsfonts}
\usepackage[rgb]{xcolor}
\usepackage{graphicx}
\def\sqr#1#2{{\vcenter{\vbox{\hrule height.#2pt
              \hbox{\vrule width.#2pt height#1pt \kern#1pt \vrule width.#2pt}
          \hrule height.#2pt}}}}

\usepackage{amsmath,amssymb,amsthm}

\def\dbF{\hbox{\rm l\negthinspace F}}

\def\dbP{\hbox{\rm l\negthinspace P}}

\def\dbR{\hbox{\rm l\negthinspace R}}

\def\a{\alpha}
\def\b{\beta}

\def\th{\theta}

\def\G{\Gamma}

\def\L{\Lambda}

\def\cF{{\cal F}}

\def\cU{{\cal U}}

\def\sqr#1#2{{\vcenter{\vbox{\hrule height.#2pt
              \hbox{\vrule width.#2pt height#1pt \kern#1pt \vrule width.#2pt}
              \hrule height.#2pt}}}}

\def\dbR{{\mathop{\rm l\negthinspace R}}}

\def\3n{\negthinspace \negthinspace \negthinspace }
\def\2n{\negthinspace \negthinspace }
\def\1n{\negthinspace }

\def\dbF{{\mathop{\rm l\negthinspace F}}}

\def\dbP{{\mathop{\rm l\negthinspace P}}}
\def\dbR{{\mathop{\rm l\negthinspace R}}}
\def\={\buildrel \triangle \over =}

\def\a{\alpha}
\def\b{\beta}

\def\th{\theta}

\def\G{\Gamma}

\def\L{\Lambda}

\def\cF{{\cal F}}

\def\cU{{\cal U}}

\def\no{\noindent}

\def\ms{\medskip}

\def\exp{\mathop{\rm exp}}
\def\sup{\mathop{\rm sup}}

\def\inf{\hbox{\rm inf$\,$}}

\def\as{\hbox{\rm a.s.{ }}}

\def\|{\Big |}
\def\({\Big (}
\def\){\Big )}
\def\[{\Big[}
\def\]{\Big]}
\def\be{\begin{equation}}
\def\bel{\begin{equation}\label}
\def\ee{\end{equation}}
\def\bt{\begin{theorem}}
\def\bcd{\begin{condition}}
\def\ecd{\end{condition}}
\def\et{\end{theorem}}
\def\bc{\begin{corollary}}
\def\ec{\end{corollary}}
\def\bde{\begin{definition}}
\def\ede{\end{definition}}
\def\bl{\begin{lemma}}
\def\el{\end{lemma}}
\def\bp{\begin{proposition}}
\def\ep{\end{proposition}}
\def\br{\begin{remark}}
\def\er{\end{remark}}
\def\ba{\begin{array}}
\def\ea{\end{array}}
\def\ed{\end{document}}

\def\square#1{\vbox{\hrule\hbox{\vrule height#1%
     \kern#1\vrule}\hrule}}
\def\rectangle#1#2{\vbox{\hrule\hbox{\vrule height#1%
     \kern#2\vrule}\hrule}}
\def\qed{\hfill \vrule height7pt width3pt depth0pt}

\font\tenbb=msbm10 \font\sevenbb=msbm7 \font\fivebb=msbm5

\newfam\bbfam
\scriptscriptfont\bbfam=\fivebb \textfont\bbfam=\tenbb
\scriptfont\bbfam=\sevenbb

\newtheorem{lemma}{Lemma}
\newtheorem{remark}{Remark}

\newtheorem{theorem}{Theorem}
\newtheorem{corollary}{Corollary}

\newtheorem{definition}{Definition}
\newtheorem{proposition}{Proposition}
\newtheorem{condition}{Assumption}

\begin{document}

\title{ Risk-Sensitive Mean-Field Type Control under Partial Observation}

\author{Boualem Djehiche\thanks{KTH Royal Institute of Technology, 
\emph{E-mail: boualem@math.kth.se} } \and Hamidou Tembine\thanks{New York University, 
 \emph{E-mail: tembine@nyu.edu} }}

\maketitle
\thispagestyle{empty}
\pagestyle{empty}
\begin{abstract}
We establish a stochastic maximum principle (SMP) for control problems of partially observed diffusions of mean-field type with risk-sensitive performance functionals.
\end{abstract}


\begin{center}
\begin{minipage}{14.6cm}
\end{minipage}
\end{center}
\begin{center}
\begin{minipage}{14.6cm}
\end{minipage}
\end{center}

\section{Introduction}

In optimal control problems for diffusions of mean-field type the performance functional, drift and diffusion coefficient depend not only on the state and the control but also on the probability distribution of state-control pair. The mean-field coupling makes the control problem time-inconsistent in the sense that the Bellman Principle is no longer valid, which motivates the use of the stochastic maximum (SMP) approach to solve this type of optimal control problems instead of trying extensions of the dynamic programming principle (DPP). This class of control problems has been studied by many authors including \cite{b1,b2,b3, carmona, hosking, li}. The performance functionals considered in these papers have been of risk-neutral type i.e. the running cost/profit terms are expected values of stage-additive payoff functions. Not all behavior, however, can be captured by risk-neutral performance. One way of capturing risk-averse and risk-seeking behaviors is by exponentiating the performance functional before expectation (see \cite{Jac73}). 

\ms The first paper that we are aware of and which deals with risk-sensitive optimal control in a mean field context is \cite{tembine2014}. Therein, the authors derive a verification theorem for a risk-sensitive mean-field game whose underlying dynamics is a Markov diffusion, using a matching argument between a system of Hamilton-Jacobi-Bellman (HJB) equations and the Fokker-Planck equation. This matching argument freezes the mean-field coupling in the dynamics, which yields a standard risk-sensitive HJB equation for the value-function. The mean-field coupling is then retrieved through the Fokker-Planck equation satisfied by the marginal law of the optimal state. 

\ms In a recent paper \cite{DT-1}, the authors have established a risk-sensitive SMP for mean-field type control. The risk-sensitive control problem was first reformulated in terms of an augmented state process and terminal payoff  problem. An intermediate stochastic maximum principle was then obtained by applying the SMP of (\cite{b3}, Theorem 2.1.) for loss functionals without running cost but  with augmented state in higher dimension and complete observation of the state. Then, the intermediate first- and second-order adjoint processes are  transformed to a more simpler form using a logarithmic transformation derived in \cite{Karoui-Ham}.

\ms Optimal control of partially observed diffusions (without mean-field coupling) has been studied by many authors including the non-exhaustive references \cite{ fleming, Davis, kwakernaak, bensoussan, haussmann, baras-elliot-kohlmann,  whittle, zhou,  li-tang, charala, charala2, tang,  Huang}, using both the DPP and SMP approaches. \cite{tang} derives an SMP for the most general model of optimal control of partially observed diffusions under risk-neutral performance functionals. Recently, Wang {\it et al.} \cite{Wang}, extended the SMP for  partially observable optimal control of diffusions for risk-neutral performance functionals of mean-field type. 

\ms The purpose of this paper is to establish a stochastic maximum principle for a class of risk-sensitive mean-field type control problems under partial observation. Following the above mentioned papers of optimal control under partial observation, in particular \cite{tang}, our strategy is to transform the partially observable control problem into a completely observable one and then apply the approach suggested in \cite{DT-1} to derive the suitable the risk-sensitive SMP. To the best to our knowledge, the risk-sensitive maximum principle under partial observation without passing through the DPP, and in particular, for mean-field type controls has not been established in earlier work.

\ms The paper is organized as follows. In Section \ref{sec:results}, we present the model and state the partially observable risk-sensitive SMP which constitutes the main result, whose prove is displayed in Section \ref{sec: main result}. Finally, in Section \ref{sec:example}, we apply the risk-sensitive SMP to the linear-exponential- quadratic setup under partial observation. To streamline the presentation, we only consider the one-dimensional case. The extension to the multidimensional case is by now straightforward. Furthermore, we  consider diffusion models where the control enters only the drift coefficient, which leads to an SMP with only one pair of adjoint processes. The general Peng-type SMP  can be derived following e.g. \cite{tang} and \cite{DT-1}.

\section{Statement of the problem}\label{sec:results}
 Let $T>0$ be a fixed time horizon and $(\Omega,{\cal F},\dbF, \dbP)$  be a given filtered  probability space on which are defined two independent  standard one-dimensional 
Brownian motions $W=\{W_s\}_{s\geq0}$ and $Y=\{Y_s\}_{s\geq0}$. Let  $\mathcal{F}_t^{W}$ and  $\mathcal{F}_t^{Y}$ be the $\dbP$-completed natural filtrations generated by $W$ and $Y$, respectively. Set  $\dbF:=\{{\mathcal{F}}_s,\ 0\leq s \leq T\}$, where, $\mathcal{F}_t=\mathcal{F}_t^{W} \vee \mathcal{F}_t^{Y}$.

\ms We consider the stochastic controlled system of mean-field type with partial observation which has the following characteristics:

\ms $(i)$ An admissible control $u$ is an $\dbF^{Y}$-adapted process with values in a non-empty subset (not necessarily convex) $U$ of $\dbR$ and satisfies $E[\int_0^T|u(t)|^2dt]<\infty$. We denote the set of all admissible controls by $\mathcal{U}$. The control $u$ is called partially observable.
 
\ms $(ii)$  Given a control process $u\in \mathcal{U}$, the controlled state process $x^u(\cdot)$ can only be partially observed through $Y$, which we call the observation process, via the dynamics 
\begin{equation}\label{SDEu-y}
\left\{\begin{array}{lll}
dY_t=\b(t,x^u(t))dt+  d\widetilde{W}_{t}, \\ Y_0=0,\\
\end{array}\right.
\end{equation}
on $(\Omega,{\cal F},\dbF, \dbP)$, where $\b(t,x): [0,T] \times \dbR\longrightarrow \dbR,$ is a Borel measurable function. A more general model of the function $\b$ would be to let it depend on the control $u$ and be of mean-field type. To keep the presentation simpler, we skip these cases in this paper. But, the main results do extend to this case.

\ms $(iii)$ Under a probability measure $ \dbP^u$, the state process $x^u(\cdot)$ satisfies the following SDE of mean-field type  
\begin{equation}\label{SDEu-x}
\left\{\begin{array}{lll}
dx^u(t)=b(t,x^u(t),E^u[x^u(t)], u(t))dt+\sigma(t,x^u(t),E^u[x^u(t)])dW_{t}\\ \qquad\quad+\a(t,x^u(t),E^u[x^u(t)])d\widetilde{W}_{t}, \\ x^u(0)=x_0,
\end{array}\right.
\end{equation}
where, $W$ and $\widetilde W$ are two independent Brownian motions,  $x_0$ is assumed  real constant, and  $E^{u}$ denotes the expectation with the respect to probability measure $ \dbP^{u}$, 
\begin{equation*}
b(t,x,m,u): \,\,[0,T] \times \dbR\times\dbR\times U\longrightarrow \dbR,
\end{equation*}
and 
\begin{equation*}
\a(t,x,m),\,\,\sigma(t,x,m): \,\,[0,T] \times \dbR\times\dbR\longrightarrow \dbR.
\end{equation*}

\ms\no  The objective is to characterize admissible controls which minimize the risk-sensitive cost functional associated with (\ref{SDEu-x}) given by
\begin{equation}\label{rs-cost}
J^{\theta}(u(\cdot))=E^{u}\left[\exp{\left(\theta\left[\int_0^Tf(t,x^u(t),E^u[x^u(t)], u(t))\,dt+ h(x^u(T),E^u[x^u(T)])\right]\right)}\right],
\end{equation}
where, $\theta$ is the risk-sensitivity index, 
\begin{equation*}\begin{array}{lll}
f(t,x,m,u): \,\,[0,T] \times \dbR\times\dbR\times U\longrightarrow \dbR,\ \ h(x,m): \,\,\dbR\times\dbR\longrightarrow \dbR, \\ t\in[0, T],\ x\in \dbR,\ m\in \dbR,\ u \in U.
\end{array}
\end{equation*}

\ms\no Any $\bar u(\cdot)\in {\cal U}$ which satisfies
\be\label{rs-opt-u}
  J^{\theta}(\bar u(\cdot))=\inf_{u(\cdot)\in {\cal U}}J^{\theta}(u(\cdot))
\ee
is called a risk-sensitive optimal control under partial observation.

\ms\no  Let $\Psi_T=\int_0^T f(t, x(t),E^u[x(t)], u(t)) dt+h(x(T), E^u[x(T)])$ and consider the payoff functional given by 
$$
\widetilde\Psi_{\theta}:=\frac{1}{\theta}\log  E^u e^{\theta \Psi_T}.
$$
When the risk-sensitive index $\theta$ is small, the loss functional $\widetilde\Psi_{\theta}$ can be expanded as
$$
E^u[ \Psi_T] +\frac{\theta}{2}\mbox{var}_u(\Psi_T)+O(\theta^2),
$$
where,  $\mbox{var}_u(\Psi_T)$ denotes the variance of  $\Psi_T$ w.r.t. $ \dbP^u$. If $\theta<0$ , the variance of $\Psi_T$, as a measure of risk,  improves the performance $\widetilde\Psi_{\theta}$, in which case the optimizer is called {\it risk seeker}. But,  when $\theta>0$,  the variance of $\Psi_T$ worsens the performance $\widetilde\Psi_{\theta}$,   in which case the optimizer is called {\it risk averse}.
The risk-neutral loss functional $E^u[\Psi_{T}]$ can be seen as a limit of risk-sensitive functional $ \widetilde\Psi_{\theta}$ when $\theta\rightarrow 0$.

\ms\no Introduce the density process defined on $(\Omega,{\cal F},\dbF, \dbP)$ by
\be\label{density-exp}
\rho^u(t):=\exp{\left \{\int_0^t \beta(s,x^u(s))dY_s-\frac{1}{2}\int_0^t|\beta(s,x^u(s))|^2ds \right\}},
\ee
which solves the linear SDE
\begin{equation}\label{density}
d\rho^u(t)=\rho^u(t)\beta(t,x^u(t)) dY_t,\,\, \ \rho^u(0)=1.\
\end{equation}
Assuming  the function $\b$ bounded (see Assumption \ref{cond1}, below), $\rho$ is a uniformly integrable martingale such that, for every $k\ge 2$,
\be\label{int-density}
E[\sup_{0\le t\le T}(\rho_t^u)^k]\le C,
\ee
where, $C$ is a constant which depends only on the bound of $\beta$, $p$ and $T$.  By Girsanov's Theorem, $d\dbP^{u}=\rho^u(T)d\dbP$. Moreover, $\dbP$ and $\dbP^u$ are equivalent measures.  This relationship between $\dbP$ and $\dbP^u$ enables us to merge (\ref{SDEu-y}) with  (\ref{SDEu-x}) and obtain the controlled state process $(\rho^u,x^u)$ as a weak solution on $(\Omega,{\cal F},\dbF, \dbP)$ of the following dynamics:
\begin{equation}\label{SDEu}
\left\{\begin{array}{lll}
d\rho^u(t)=\rho^u(t)\beta(t,x^u(t)) dY_t,\\ 
dx^u(t)=\left\{b(t,x^u(t),E[\rho^u(t)x^u(t)], u(t))-\a(t,x^u(t),E[\rho^u(t)x^u(t)])\b(t,x^u(t))\right\}dt\\  \qquad\quad+\sigma(t,x^u(t),E[\rho^u(t)x^u(t)])dW_{t} +\a(t,x^u(t),E[\rho^u(t)x^u(t)])dY_{t}, \\ \rho^u(0)=1,\,\, x^u(0)=x_0.
\end{array}\right.
\end{equation}
Moreover, the associated risk-sensitive cost functional (\ref{rs-cost}) becomes 
\begin{equation}\label{cost}
J^{\theta}(u(\cdot))=E\left[\rho^u(T)e^{\theta\left[\int_0^Tf(t,x^u(t),E[\rho^u(t)x^u(t)], u(t))\,dt+ h(x^u(T),E[\rho^u(T)x^u(T)])\right]}\right].
\end{equation}
We have recast the partially observable control problem into a completely observable control problem of the state process $(\rho^u,x^u)$ which for instance boils down to characterizing the controls $\bar u(\cdot)\in \mathcal{U}$ which satisfy (\ref{rs-opt-u}), where the cost functional $J^{\th}$ is given by (\ref{cost}), subject to the dynamics  $(\rho^u, x^u)$ solution of (\ref{SDEu}).

\ms\no The main result of this paper is a stochastic maximum principle (SMP) in terms of necessary optimality conditions for the problem (\ref{rs-opt-u}), subject to (\ref{SDEu})-(\ref{cost}).

\ms We will make the following assumption.
\bcd\label{cond1} The functions $b, \sigma,\a,\b, f, h$ are twice continuously differentiable with respect to $(x,m)$.
Moreover, these functions and their first derivatives with respect to $(x,m)$ are continuous in $(x,m, u)$, and bounded.
\ecd

\no To keep the presentation less technical, we impose these assumptions although they are restrictive and can be made weaker.

\no  Under these assumptions, in view of Girsanov's theorem and  \cite{jourdain}, Proposition 1.2.,  for each $u\in {\cal U}$, the SDE (\ref{SDEu}) admits a unique weak solution $(\rho^u, x^u)$.

\ms\no We now state  an SMP to characterize optimal controls $\bar u(\cdot)\in \mathcal{U}$ which minimize (\ref{cost}), subject to (\ref{SDEu}). Let $(\bar\rho,\bar x):=(\rho^{\bar u},x^{\bar u})$ denote the corresponding state process, solution of (\ref{SDEu}).

\ms\no We introduce the following notation.
\be\label{notation-2}\begin{array}{lll}
X:=\left(\begin{array}{lll} \rho\\ x \end{array}\right),\,\,\, \bar X:=\left(\begin{array}{lll} \bar\rho\\ \bar x\\ \end{array}\right),\,\,\, \phi(X)=\rho x,\,\,\, \phi(\bar X)=\bar\rho \bar x,\,\,\,  X_0=\bar X_0:=\left(\begin{array}{lll} 1\\ x_0\end{array}\right),\,\,\, B_t:=\left(\begin{array}{lll} Y_t \\ W_t\end{array}\right),\\ c(t,x,m,u):= b(t,x,m,u)-\a(t,x,m)\b(t,x),\\ 
F(t,X,m,u):=\left(\begin{array}{lll} 0\\ c(t,x,m,u)\end{array}\right),\,\,\,
G(t,X,m):=\left(\begin{array}{lll} \rho\b(t,x) & 0\\ \a(t,x,m) & \sigma(t,x,m)\end{array}\right),
\end{array}
\ee

\ms\no We define the  risk-neutral Hamiltonian  as follows. For $(p,q)\in \dbR^2\times\dbR^{2\times 2}$, 
\be\label{rn-hamiltonian-1}
H(t,X,m,p,q,u):=\langle F(t,X,m,u),p\rangle+\mbox{tr}(G^*(t,X,m)q)-f (t, x,m,u), 
\ee
where, $'*'$ denotes the transposition operation of a matrix or a vector.

\ms\no We also introduce the risk-sensitive Hamiltonian: for $\theta \in \dbR$ and $(p,q, \ell)\in \dbR^2\times \dbR^{2\times 2}\times \dbR^2$,
\begin{equation}\label{rs-hamiltonian-1}\begin{array}{lll}
H^{\theta}(t, X,m, u, p, q,\ell):= \langle F(t,X,m,u),p\rangle -f (t, x,m, u)+\mbox{tr}(G^*(t,X,m)(q +\theta\ell p^*)).
\end{array}
\end{equation}
We have $H=H^0$.

\ms\no Setting 
\be\label{pq-notation-1}
\ell:=\left(\begin{array}{lll} \ell_1\\ \ell_2 \end{array}\right), \quad p:=\left(\begin{array}{lll} p_1\\ p_2 \end{array}\right),\quad q:=\left(\begin{array}{lll}  q_{11} & q_{12}\\   q_{21} & q_{22} \end{array}\right),
\ee
the explicit form of the Hamiltonian (\ref{rs-hamiltonian-1}) reads
\be\label{rs-hamiltonian-1-1}\begin{array}{lll}
H^{\theta}(t, X,m, u, p, q,\ell):=c(t,x,m,u)p_2-f(t,x,m,u)+\rho\b(t,x)(q_{11}+\theta \ell_1 p_1)\\ \qquad\qquad\qquad\qquad\quad +\a(t,x,m)(q_{21}+\theta\ell_2 p_1)+\sigma(t,x,m)(q_{22}+\theta\ell_2 p_2).
\end{array}
\ee
Setting $\theta=0$ in (\ref{rs-hamiltonian-1-1}), we obtain the explicit form of the Hamiltonian (\ref{rn-hamiltonian-1}):
\be\label{rn-hamiltonian-1-1}\begin{array}{lll}
H(t, X,m, u, p, q):=c(t,x,m,u)p_2-f(t,x,m,u)+\rho\b(t,x)q_{11}\\ \qquad\qquad\qquad\qquad\quad +\a(t,x,m)q_{21}+\sigma(t,x,m)q_{22}.
\end{array}
\ee
With the obvious notation for the derivatives of the functions $b,\a,\b, \sigma,f,h$, w.r.t. the arguments $x$ and $m$, we further set
\be\label{dH-rn-1}\left\{\begin{array}{lll}
H^{\theta}_x(t, X,m, u, p, q):=c_x(t,x,m,u)p_2-f_x(t,x,m,u)+\rho\b_x(t,x)(q_{11}+\theta \ell_1 p_1) \\ \qquad\qquad\qquad\qquad\quad+\a_x(t,x,m)(q_{21}+\theta\ell_2 p_1)+\sigma_x(t,x,m)(q_{22}+\theta\ell_2 p_2),\\ 
H^{\theta}_m(t, X,m, u, p, q):=c_m(t,x,m,u)p_2-f_m(t,x,m,u)\\ \qquad\qquad\qquad\qquad\quad+\a_m(t,x,m)(q_{21}+\theta\ell_2 p_1)+\sigma_m(t,x,m)(q_{22}+\theta\ell_2 p_2),\\
H^{\theta}_{\rho}(t,X,m,u,p,q)=\b(t,x)(q_{11}+\theta \ell_1 p_1).
\end{array}
\right.
\ee


\ms\no With this notation, the system (\ref{SDEu}) can be rewritten in the following compact form
\be\label{SDEu-3}\left\{ 
\begin{array}{lll}
dX(t)=F (t, X(t),E[\phi(X(t))],u(t)) dt + G(t, X(t),E[\phi(X(t))]) dB_t,\\ X(0)=X_0,
\end{array}
\right.
\ee

\ms\no We define the risk-neutral Hamiltonian associated with random variables $X$ such that $\phi(X)\in L^1(\Omega,\cF,\dbP)$ as follows (with the obvious abuse of notation): For $(p,q)\in \dbR^2\times\dbR^{2\times 2}$, 
\be\label{rn-hamiltonian}
H(t,X,p,q,u):=\langle F(t,X,E[\phi(X)],u),p\rangle -f(t,x,E[\phi(X)],u)+\mbox{tr}(G^*(t,X,E[\phi(X)])q), 
\ee

\ms\no We also introduce the risk-sensitive Hamiltonian: for $\theta \in \dbR$ and $(p,q, \ell)\in \dbR^2\times \dbR^{2\times 2}\times \dbR^2$,
\begin{equation}\label{rs-hamiltonian}\begin{array}{lll}
H^{\theta}(t, \rho, x, u, p, q,\ell)=H^{\theta}(t, X, u, p, q,\ell):= \langle F(t,X,E[\phi(X)],u),p\rangle -f (t,x ,E[\phi(X)], u)\\ \qquad\qquad\qquad\qquad+\mbox{tr}(G^*(t,X,E[\phi(X)])(q +\theta\ell p^*)).
\end{array}
\end{equation}

\no\ms 
For $\phi\in \{b,c, \sigma,\a,\b, f, h\}$ and $u\in U$, we introduce the stochastic processes
\be\label{notation-3}\begin{array}{llll}
\phi_x(t):=\phi_x(t,\bar x(t),E[\bar\rho(t)\bar x(t)],\bar u(t)),\,\,\,
\phi_m(t):=\phi_m(t,\bar x(t),E[\bar\rho(t)\bar x(t)],\bar u(t)).
\end{array}
\ee

\ms\no Let 
\begin{equation}\label{psi}
\psi^{\theta}_T:=\bar\rho(T)\exp{\theta \left[\int_0^Tf(t,\bar x(t), E[\bar\rho(t)\bar x(t)], \bar u(t)) dt+h(\bar x(T), E[\bar\rho(T)\bar x(T)])\right]}.
\end{equation}

\ms\no We introduce the adjoint equations involved in the risk-sensitive SMP for our control problem.
\be\label{rs-pq}
\left\{\begin{array}{lll}
d\hat p(t)=-\left(\begin{array}{ccc} H^{\theta}_{\rho}(t)+\frac{\bar x(t)}{v^{\theta}(t)}E[v^{\theta}(t)H^{\theta}_m(t)] \\ H^{\theta}_x(t)+\frac{\bar \rho(t)}{v^{\theta}(t)}E[v^{\theta}(t)H^{\theta}_m(t)]\end{array}\right)dt+\hat q(t)(-\theta \ell(t)dt+dB_t),\\
dv^{\theta}(t)=\theta v^{\theta}(t)\langle \ell(t),dB_t\rangle,\\
\hat p(T)=-\left(\begin{array}{ccc} (\theta\bar\rho(T))^{-1}\\ h_x(T) \end{array}\right)-\left(\begin{array}{ccc} \bar x(T)\\ \bar \rho(T)\end{array}\right)\frac{1}{\psi^{\theta}_T}E[\psi^{\theta}_Th_m(T)],\\ v^{\theta}(T)=\psi^{\theta}_T,
\end{array}
\right.
\ee
\ms\no  where, in view of (\ref{dH-rn-1}) and (\ref{notation-3}), for $k=\rho,x,m$,
\be\label{rs-dH}\begin{array}{lll}
H^{\theta}_k(t):=\langle F_k(t,\bar X(t),E[\phi(\bar X(t))],\bar u(t)), \hat p(t)\rangle-f_k (t, \bar x(t),E[\phi(\bar X(t))], \bar u(t))\\  \qquad\qquad\qquad\qquad+\mbox{tr}(G_k^*(t,\bar X(t),E[\phi(\bar X(t))])(\hat q(t) +\theta\ell {\hat p}^*(t)).
\end{array}
\ee
We note that the processes $(\hat p,\hat q,\ell)$ may depend on the sensitivity index $\theta$. To ease notation, we omit to make this dependence explicit.

\ms\no Below, we will show that, under Assumption 1, (\ref{rs-pq}) admits a unique $\dbF$-adapted solution $(\hat p,\hat q,v^{\theta},\ell)$ such that
\be\label{rs-pq-bounds}
E\left[\sup_{t\in[0,T]}|\hat p(t)|^2+\sup_{t\in[0,T]}|v^{\theta}(t)|^2+\int_0^T \left(|\hat q(t)|^2+|\ell(t)|^2\right) dt\right]<\infty.
\ee 

\ms\no Moreover, 
\begin{lemma}\label{L-mart} The process defined on $(\Omega,{\cal F},\dbF, \dbP)$  by 
\be\label{v-L}
L^{\theta}_t:=\frac{v^{\theta}(t)}{v^{\theta}(0)}=\exp{\left(\int_0^t \theta\langle \ell(s),dB_s\rangle - \frac{\theta^2}{2}\int_0^t |\ell (s)|^2 ds\right)}, \quad 0\le t\le T.
\ee
is a uniformly integrable $\dbF$-martingale. 
\end{lemma}

\no The process $L^{\theta}$ defines a new probability measure $\dbP^{\theta}$ equivalent to $\dbP$ by setting $L^{\theta}_t:=\frac{d\dbP^{\theta}}{d\dbP}\|_{{\cal F}_t}$. By Girsanov's theorem, the process $B_t^{\theta}:=B_t-\theta\int_0^t \ell(s)ds,\,\, 0\le t\le T$ is a $\dbP^{\theta}$-Brownian motion.

\ms\no The following theorem is the main result of the paper.

\begin{theorem}\label{main result}$(${\bf Risk-sensitive maximum principle}$)$
Let Assumption \ref{cond1} hold. If  $(\bar \rho(\cdot),\bar x(\cdot),\bar  u(\cdot))$  is an optimal solution of the risk-sensitive control problem (\ref{rs-opt-u})-(\ref{SDEu}), then there are two pairs of $\dbF$-adapted processes $(v^{\theta}, \ell)$ and $(\hat p,\hat q)$ which satisfy (\ref{rs-pq})-(\ref{rs-pq-bounds}), such that
\be\label{rs-VI}\begin{array}{lll}
E^{\theta}[H^{\theta}(t, \bar \rho(t),\bar x(t), \hat p(t), \hat q(t),\ell(t),u)-H^{\theta}(t, \bar\rho(t),\bar x(t), \hat p(t), \hat q(t),\ell(t),\bar u(t))|\mathcal{F}^Y_t] \leq 0,
\end{array}
\ee
 for all $u\in U,$ almost every $t$ and $\dbP^{\theta}-$almost surely.

\ms\no  Here, $E^{\theta}[\,\cdot\,]$ denotes the expectation w.r.t. $\dbP^{\theta}$.
\end{theorem}

\begin{remark}\label{R-main}
\ms\no The boundedness assumption of the involved coefficients and their derivatives imposed in Assumption 1, in Theorem \ref{main result}, guarantees the solvability of the system of forward-backward SDEs (\ref{SDEu})-(\ref{rs-pq}). In fact Theorem \ref{main result} applies provided we can solve system of forward-backward SDEs (\ref{SDEu}) and (\ref{rs-pq}). A typical example of such a situation is the classical Linear-Quadratic (LQ) control problem (see Section \ref{sec:example} below), in which the involved coefficients are at most quadratic, but not necessarily bounded.  

\end{remark}

\section{Proof of the main result}\label{sec: main result}
In this section we give a proof of Theorem \ref{main result} displayed in several steps.
\subsection{An intermediate SMP for mean-field type control}\label{sec:intermediate}
 In this subsection we first reformulate the risk-sensitive control problem (\ref{rs-opt-u})-(\ref{SDEu}) in terms of an augmented state process and terminal payoff  problem. An intermediate stochastic maximum principle is then obtained by applying the SMP obtained in (\cite{b1}, Theorem 3.1 and \cite{b3}, Theorem 2.1) for loss functionals without running cost. Then,
we transform the intermediate first-order adjoint processes to a more simpler form.  The  mean-field type control  problem  (\ref{cost}) under the dynamics (\ref{SDEu}) is equivalent to 

\be\label{smp2}
\inf_{u(\cdot)\in \cU}  E\left[\rho(T) e^{\theta \left[h(x(T), E[\rho(T)x(T)])+ \xi(T)\right]}\right],
\ee
subject to
\begin{eqnarray}
\label{SDEu-2}
\left\{
\begin{array}{lll} 
d\rho(t)=\rho(t)\beta(t,x(t)) dY_t,\\ 
dx(t)=\left\{b(t,x(t),E[\rho(t)x(t)], u(t))-\a(t,x(t),E[\rho(t)x(t)])\b(t,x(t))\right\}dt\\  \qquad\quad+\sigma(t,x(t),E[\rho(t)x(t)])dW_{t} +\a(t,x(t),E[\rho(t)x(t)])dY_{t},\\
 d\xi= f(t, x(t),E[\rho(t) x(t)], u(t))dt, \\ \rho(0)=1,\,\, x(0)=x_0,\,\xi(0)=0,
 \end{array}
\right.
\end{eqnarray}

\ms\no We introduce the following notation.
\be\label{notation-4}\begin{array}{ccc}
R:=\left(\begin{array}{ccc} \rho\\ x \\ \xi \end{array}\right)=\left(\begin{array}{ccc} X \\ \xi \end{array}\right),\,\, \bar R:=\left(\begin{array}{ccc} \bar\rho\\ \bar x\\ \bar\xi \end{array}\right)=\left(\begin{array}{ccc} \bar X \\ \bar \xi \end{array}\right),\,\,\,  R_0=\bar R_0:=\left(\begin{array}{ccc} X_0,\\ 0 \end{array}\right),\\
\L(t,R,m,u):=\left(\begin{array}{ccc} F(t,X,m,u)\\ f(t,x,m,u)\end{array}\right),\,\,\,
\G(t,R,m):=\left(\begin{array}{ccc} G(t,X,m,u)\\ 0  \end{array}\right), \\ \phi(R)=\phi(X),\,\,\, \phi(\bar R)=\phi(\bar X).
\end{array}
\ee

\ms\no With this notation, the system (\ref{SDEu-2})  can be rewritten in the following compact form
\be\label{SDEu-3}\left\{ 
\begin{array}{lll}
dR(t)=\L (t, R(t),E[\phi(R(t))],u(t)) dt + \G(t, R(t),E[\phi(R(t))]) dB_t,\\ R(0)=R_0,
\end{array}
\right.
\ee 
and the  risk-sensitive cost functional  (\ref{cost}) is given by
\begin{equation}\label{cost-3}
J^{\theta}(u(\cdot)):=E[\Phi\left(R(T),E[\phi(R(T))]\right)],
\end{equation}
where 
\be\label{cost-f}
\Phi\left(R(T),E[\phi(R(T))]\right):=\rho(T) \exp\left(\theta h(x(T), E[\rho(T)x(T)]) +\theta \xi(T)\right).
\ee
We define the Hamiltonian associated with random variables $R$ such that $\phi(R)\in L^1(\Omega,\cF,\dbP)$ as follows. For $(p,q)\in \dbR^3\times\dbR^{3\times 3}$ 
\be\label{H-f}
H^e(t,R,p,q,u):=\langle \L(t,R,E[\phi(R)],u),p\rangle+\mbox{tr}(\G^*(t,R,E[\phi(R)])q), 
\ee
where, $\G^*$ denotes the transpose of the matrix $\G$.

\ms\no Setting 
\be\label{pq-notation}
p:=\left(\begin{array}{lll} p_1\\ p_2 \\ p_3 \end{array}\right),\quad q:=\left(\begin{array}{lll}  q_{11} & q_{12}\\   q_{21} & q_{22}\\ q_{31} & q_{32} \end{array}\right),
\ee
the explicit form of the Hamiltonian (\ref{H-f}) reads
\be\label{H-f-1}\begin{array}{lll}
H^e(t,\rho,x,\xi,p,q,u):=H^e(t,R,p,q,u)=c(t,x,E[\rho x],u)p_2+f(t,x,E[\rho x],u)p_3\\ \qquad\qquad\qquad \quad+\sigma(t,x,E[\rho x])q_{22}+\rho\b(t,x) q_{11}+\a(t,x,E[\rho x]) q_{21}.
\end{array}
\ee

\ms\no  In view of (\ref{notation-3}), we set, for $u\in U$,
\be\label{dH-f}\begin{array}{lll}
H^e_x(t):=\langle \L_x(t,\bar X(t),E[\phi(\bar X(t))],\bar u(t)), p\rangle+\mbox{tr}(\G_x^*(t,\bar X(t),E[\phi(\bar X(t))])q),
\end{array}
\ee
for $k=\rho,x,m$,
\ms\no whose explicit form is 
\be\label{dH-f-1}\left\{\begin{array}{lll}
H^e_x(t):=c_x(t)p_2(t)+f_x(t)p_3(t)+\sigma_x(t)q_{22}(t)+\bar\rho(t)\b_x(t) q_{11}(t)+\a_x(t)q_{21}(t),\\ 
H^e_m(t):=c_m(t)p_2(t)+f_m(t)p_3(t)+\sigma_m(t)q_{22}(t)+\bar\rho(t)\b_m(t) q_{11}(t)+\a_m(t) q_{21}(t),\\
H^e_{\rho}(t)=\b(t,\bar x(t))q_{11}(t).
\end{array}
\right.
\ee

\ms\no We may apply the  SMP for risk-neutral mean-field type control (cf. \cite{b1}, Theorem 3.1 and \cite{b3}, Theorem 2.1) to the augmented state dynamics $(\rho, x,\xi)$ to derive the first order adjoint equation:
\be\label{pq-f}
\left\{\begin{array}{lll}
dp(t)=-\left(\begin{array}{ccc} H^e_{\rho}(t)+\bar x(t)E[H^e_m(t)]\\ H^e_x(t)+\bar \rho(t)E[H^e_m(t)]\\ 0\end{array}\right)dt+q(t)dB_t,\\
p(T)=-\theta \psi^{\theta}_T\left(\begin{array}{ccc} (\theta\bar\rho(T))^{-1}\\ h_x(T)\\ 1\end{array}\right)-\theta\left(\begin{array}{lll} \bar x(T)\\ \bar \rho(T)\\ 0\end{array}\right)E[\psi^{\theta}_Th_m(T)].
\end{array}\right.
\ee
This is a system of linear backward SDEs with mean-field type which, in view of (\cite{Buckdahn1}, Theorem 3.1), under Assumption \ref{cond1}, admits a unique  $\dbF$-adapted solution $(p,q)$ which satisfies  
\be\label{pq-f-bound}
E\left[\sup_{t\in[0,T]}|p(t)|^2+\int_0^T |q(t)|^2 dt\right]<\infty.
\ee 
where, $|\cdot|$ denotes the usual Euclidean norm with appropriate dimension.

\ms\no We may apply the SMP for SDEs of mean-field type control  from  (\cite{b1}, Theorem 3.1 and \cite{b3}, Theorem 2.1)  together with the SMP for risk-neutral partially observable SDEs derived in (\cite{tang}, Theorem 2.1) to obtain the following SMP.
\begin{proposition}\label{rn-MSP} 
Let Assumption \ref{cond1} hold. If $(\bar R(\cdot), \bar  u(\cdot))$ is an optimal solution of the risk-neutral control problem  (\ref{smp2}) subject to the dynamics (\ref{SDEu-2}),  then there is a unique pair of $\dbF$-adapted processes $(p,q)$ which  satisfies (\ref{pq-f})-(\ref{pq-f-bound})  such that
\begin{equation}\label{rn-VI}\begin{array}{lll}
E[H^e(t, \bar R(t), p(t), q(t),u)-H^e(t, \bar R(t), p(t), q(t),\bar u(t))|\mathcal{F}^Y_t] \leq 0, 
\end{array}
\end{equation} for all $u\in U,$ almost every $t$ and $\dbP-$almost surely.
\end{proposition}

\subsection{Transformation of the first order adjoint process}\label{sec:mainresult}
\no Although the result of Proposition \ref{rn-MSP} is a good SMP for the risk-sensitive mean-field type control with partial observations, augmenting the state process with the third component $\xi$  yields a system of three  adjoint equations that appears  complicated to solve in concrete situations. In this section we apply the transformation of the adjoint processes $(p,q)$ introduced in \cite{DT-1} in  such a way to get rid of the third component $(p_3,q_{31},q_{32})$ in (\ref{pq-f}) and express the SMP in terms of only two adjoint process that we denote $(\hat p,\hat q)$, where
\be\label{pq-notation-1}
\hat p:=\left(\begin{array}{ccc} \hat p_1\\ \hat p_2\end{array}\right),\quad \hat q:=\left(\begin{array}{ccc} \hat q_1\\ \hat q_2\end{array}\right),\quad \hat q_i:=(\hat q_{i1},\hat q_{i2}),\,\,\, i=1,2.
\ee
\medskip\noindent  Indeed,  noting that  from (\ref{pq-f}), we have 
$dp_3(t) =\langle q_3(t),dB_t\rangle$ and $ p_3(T)=- \theta \psi^{\theta}_T,$  the explicit solution of this backward SDE is  
\be\label{p3}
p_3(t)=- \theta E[\psi^{\theta}_T \ | \ {\cal F}_t ]=-\theta v^{\theta}(t) ,
\ee
where, 
\be\label{mg} 
v^{\theta}(t):= E[\psi^{\theta}_T \ | \ {\cal F}_t ],  \qquad 0\le t\le T. 
\ee
In particular, we have $v^{\theta}(0)= E[\psi^{\theta}_T].$ 
Therefore, in view of (\ref{p3}), it would be  natural to choose a transformation of $(p,q)$ into an adjoint process $(\hat{p}, \hat{q})$ ,
where,
$$
\hat p:=\left(\begin{array}{ccc} \hat p_1\\ \hat p_2 \\ \hat p_3\end{array}\right),\quad \hat q:=\left(\begin{array}{ccc} \hat q_{11} & \hat q_{12}\\ \hat q_{21} & \hat q_{22} \\ \hat q_{31} & \hat q_{32}\end{array}\right), 
$$
 such that
\be\label{t-p2}
\hat p_3(t)=\frac{p_3(t)}{ \theta v^{\theta}(t)}=-1, \qquad 0\le t\le T,
\ee
which would  imply that, for almost  every $0\le t\le T$, 
\be\label{t-q-3}
\hat q_3(t)=(\hat q_{31}(t), \hat q_{32}(t))=0, \,\,\, \dbP-\as
\ee which in turn reduces the number of adjoint process to those of the form given by (\ref{pq-notation-1}).

\ms\no  We consider the following transform:  
\be\label{transform}
\hat p(t):=\frac{1}{\theta v^{\theta}(t)} p(t),\qquad  0\le t\le T.
\ee
In view of (\ref{pq-f}), we have
\be\label{t-p-T}
\hat{p}(T)=-\left(\begin{array}{ccc} (\theta\bar\rho(T))^{-1}\\ h_x(T)\\ 1\end{array}\right)-\left(\begin{array}{ccc} \bar x(T)\\ \bar \rho(T)\\ 0\end{array}\right)\frac{1}{\psi^{\theta}_T}E[\psi^{\theta}_Th_m(T)].
\ee
\ms\no We should identify the processes $\hat \a$ and $\hat q$ such that
\be\label{t-p}
d\hat p(t)=-\hat\a(t)dt+\hat q(t) dB_t,
\ee
for which (\ref{t-p2}) and (\ref{t-q-3}) are satisfied.

\ms\no In order to investigate the properties of these new processes $(\hat p,\hat q)$, the following properties of the generic martingale $v^{\theta}$, used in \cite{DT-1}, are essential. We reproduce them here for the sake of completeness. Since, by Assumption \ref{cond1}, $f$ and $h$ are bounded by some constant $C>0$, we have 
\be\label{v-bound}
0< e^{-(1 + T)C\theta}\rho(T) \leq \psi^{\theta}_T \leq  e^{(1 + T)C\theta}\rho(T).
\ee
Therefore, $v^{\theta}$  is a uniformly integrable $\dbF$-martingale satisfying
\be\label{b-mg}
0< e^{-(1 + T)C\theta}\rho(t) \leq v^{\theta}(t) \leq e^{(1 + T)C\theta}\rho(t), \,\qquad  0\le t\le T.
\ee
Hence, in view of (\ref{int-density}), we have
\be\label{int-v}
E[\sup_{0\le t\le T}|v^{\theta}(t)|^2]\le C.
\ee
Furthermore, the martingale $v^{\theta}$ enjoys the following useful {\it  logarithmic transform}  established in (\cite{Karoui-Ham}, Proposition 3.1)

\be\label{mart-rep}
v^{\theta}(t)=\exp\left(\theta Z_t+\theta \int_0^t f(s,\bar x(s),E[\bar\rho(s)\bar x(s)],\bar u(s))ds\right), \quad 0\le t\le T, 
\ee 
and 
\be\label{Y-v}
v^{\theta}(0)=E[\psi^{\theta}_T]=\exp(\theta Z_0).
\ee
Moreover, the process $Z$ is the first component of the $\dbF$-adapted pair of processes $(Z,\ell)$  which is the unique solution to the following quadratic BSDE:
\be\label{Q-BSDE}
\left\{ \begin{array}{lll}
dZ_t=-\{f(t, \bar x(t), E[\bar\rho(s)\bar x(s)],\bar u(t))+\frac{\theta}{2}|\ell(t)|^2\}dt+\langle \ell(t), dB_t \rangle,\\  \\
Z_T=\frac{1}{\theta}\ln \bar\rho(T)+ h(\bar x_T, E[\bar\rho(T)\bar x(T)]).
\end{array}
\right.
\ee
where, $\ell(t)=(\ell_1(t),\ell_2(t))$ satisfies
\be\label{Y-ell}
E\left[\int_0^T |\ell(t)|^2dt\right]<\infty.
\ee
\medskip\noindent 
In particular, $v^{\theta}$ solves the following linear backward SDE
\be\label{v}
dv^{\theta}(t)=\theta v^{\theta}(t)\langle\ell(t),dB_t\rangle,\quad  v^{\theta}(T)=\psi^{\theta}_T.
\ee
Hence,

\ms {\bf Proof of Lemma \ref{L-mart}}.  In view of (\ref{int-v}), 
\be\label{v-L}
\frac{v^{\theta}(t)}{v^{\theta}(0)}=\exp{\left(\int_0^t \theta\langle \ell (s),dB_s\rangle - \frac{\theta^2}{2}\int_0^t |\ell (s)|^2 ds\right)}:=L^{\theta}_t, \quad 0\le t\le T.
\ee
is a uniformly integrable $\dbF$-martingale. \qed

\ms\no To identify the processes $\tilde\a$ and $\tilde q$ such that
\be\label{t-p}
d\hat p(t)=-\hat\a(t)dt+\hat q(t) dB_t,
\ee
we may  apply It\^o's formula to  the process ${p}(t)=\theta v^{\theta}\tilde{p}(t)$, use (\ref{pq-f}) and (\ref{v}) and identify the coefficients. We obtain

\be\label{pq-hat-coeff}
\left\{\begin{array}{lll}
\hat \a(t)=\frac{1}{\theta v^{\theta}(t)}\left(\begin{array}{ccc} H^e_{\rho}(t)+\bar x(t)E[H^e_m(t)]\\ H^e_x(t)+\bar \rho(t)E[H^e_m(t)]\\ 0\end{array}\right)+\theta \hat q(t)\ell(t),\\
\hat q(t)=\frac{1}{\theta v^{\theta}(t)}q(t)-\theta \hat p(t)\ell(t).\\
\end{array}
\right.
\ee
Therefore, 
\be\label{pq-hat-f}
\left\{\begin{array}{lll}
d\hat p(t)=-\frac{1}{\theta v^{\theta}(t)}\left(\begin{array}{ccc} H^e_{\rho}(t)+\bar x(t)E[H^e_m(t)]\\ H^e_x(t)+\bar \rho(t)E[H^e_m(t)]\\ 0\end{array}\right)dt+\hat q(t)dB^{\theta}_t,\\
\hat q(t)=\frac{1}{\theta v^{\theta}(t)}q(t)-\theta \hat p(t)\ell(t),\\
dv^{\theta}(t)=\theta v^{\theta}(t)\langle\ell(t),dB_t\rangle,\\
\hat p(T)=-\left(\begin{array}{ccc} (\theta\bar\rho(T))^{-1}\\ h_x(T)\\ 1\end{array}\right)-\left(\begin{array}{lll} \bar x(T)\\ \bar \rho(T)\\ 0\end{array}\right)\frac{1}{\psi^{\theta}_T}E[\psi^{\theta}_Th_m(T)],\\
v^{\theta}(T)=\psi^{\theta}_T.
\end{array}
\right.
\ee
where, $B_t^{\theta}:=B_t-\theta\int_0^t \ell(s)ds,\,\, 0\le t\le T$, which is, in view of (\ref{v-L}) and Girsanov's Theorem, a $\dbP^{\theta}$-Brownian motion, where $\frac{d\dbP^{\theta}}{d\dbP}\Big|_{{\cal F}_t}:=L^{\theta}_t$.

\ms\no In particular, 
$$
d\hat p_3(t)= \langle\hat q_3(t),-\theta \ell(t)dt+dB_t\rangle,\quad \hat p_3(T)=-1.
$$
Therefore, noting that $\hat p_3(t):=[\theta v^{\theta}(t)]^{-1}p_3(t)$ is square-integrable, we obtain $\hat p_3(t)=E^{\mathbb{P}^{\theta}}[\hat p_3(T)| {\cal F}_t]=-1$.  Thus, its quadratic variation $\int_0^T|\hat q_3(t)|^2dt=0,\,\,\mathbb{P}^{\theta}-\as$. This 
 implies that,  for almost every $0\le t\le T$, $\hat q_3(t)=0,\,\,  \mathbb{P}^{\theta}\,\,\mbox{and}\,\, \mathbb{P}-\as$

\ms\no Hence, we can drop the last components from the adjoint processes $(\hat p,\hat q)$ and only consider (keeping the same notation) 
\be\label{pq-hat-new}
\hat p:=\left(\begin{array}{ccc} \hat p_1\\ \hat p_2 \end{array}\right),\quad \hat q:=\left(\begin{array}{ccc} \hat q_{11} & \hat q_{12}\\ \hat q_{21} & \hat q_{22}\end{array}\right),
\ee
for which (\ref{pq-hat-f}) reduces to  the risk-sensitive adjoint equation: 
\be\label{pq-hat-rs}
\left\{\begin{array}{lll}
d\hat p(t)=-\frac{1}{\theta v^{\theta}(t)}\left(\begin{array}{ccc} H^e_{\rho}(t)+\bar x(t)E[H^e_m(t)] \\ H^e_x(t)+\bar \rho(t)E[H^e_m(t)]\end{array}\right)dt+\hat q(t)dB^{\theta}_t,\\
\hat q(t)=\frac{1}{\theta v^{\theta}(t)}q(t)-\theta \hat p(t)\ell(t),\\
dv^{\theta}(t)=\theta v^{\theta}(t)\langle\ell(t),dB_t\rangle,\\
\hat p(T)=-\left(\begin{array}{ccc} (\theta\bar\rho(T))^{-1}\\ h_x(T) \end{array}\right)-\left(\begin{array}{ccc} \bar x(T)\\ \bar \rho(T)\end{array}\right)\frac{1}{\psi^{\theta}_T}E[\psi^{\theta}_Th_m(T)],\\
v^{\theta}(T)=\psi^{\theta}_T.
\end{array}
\right.
\ee

\ms\no In view of the uniqueness of $\dbF$-adapted pairs $(p,q)$, solution of (\ref{pq-f}) and the pair $(v^{\theta},\ell)$ obtained by (\ref{Q-BSDE}) and (\ref{v}), the solution  of the system of backward SDEs (\ref{pq-hat-rs}) is unique and satisfies (\ref{rs-pq-bounds}).

\subsection{ Risk-sensitive stochastic maximum principle}
\no We may use the transform  (\ref{transform}) and (\ref{pq-hat-coeff}) to obtain the explicit form (\ref{rs-hamiltonian}) of the risk-sensitive Hamiltonian $H^{\theta}$ defined by
\be\label{rs-H-1}
H^{\theta}(t,\bar X(t),\hat p(t),\hat q(t),\ell(t), u):=\frac{1}{\theta v^{\theta}(t)}H^e(t,\bar R(t),p(t), q(t),u).
\ee
Let 
$$
\delta H^{e}(t):=H^e(t,\bar R(t), p(t), q(t),u)-H^e(t, \bar R(t), p(t), q(t),\bar u(t))
$$
and 
$$
\delta H^{\theta}(t)=H^{\theta}(t,\bar X(t),\hat p(t),\hat q(t),\ell(t), u)-H^{\theta}(t,\bar X(t),\hat p(t),\hat q(t),\ell(t), \bar u(t)).
$$
We have 
$$
E[\delta H^{e}(t)|\mathcal{F}^Y_t]=
\theta E[v^{\theta}(t)\delta H^{\theta}(t)|\mathcal{F}^Y_t]=\theta v^{\theta}(0)E^{\theta}[\delta H^{\theta}(t)|\mathcal{F}^Y_t],
$$
where, we recall that  $v^{\theta}(t)/v^{\theta}(0)=L^{\theta}_t=d\dbP^{\theta}/d\dbP|_{{\cal F}_t}$.

\ms\no
Now, since $\theta>0$  and $v^{\theta}(0)=E[\psi_T^{\theta}]>0$, the variational inequality  (\ref{rn-VI}) translates into 
\be\label{rs-VI-1}\begin{array}{lll}
E^{\theta}[H^{\theta}(t, \bar \rho(t),\bar x(t), \hat p(t), \hat q(t),\ell(t),u)-H^{\theta}(t, \bar \rho(t),\bar x(t), \hat p(t), \hat q(t),\ell(t),\bar u(t))|\mathcal{F}^Y_t] \leq 0. 
\end{array}
\ee
 for all $u\in U,$ almost every $t$ and $\dbP^{\theta}-$almost surely.  This finishes the proof of Theorem \ref{main result}.   $\qed$

\section{Illustrative Example:  Linear-quadratic risk-sensitive model under partial observation}  \label{sec:example}
To illustrate our approach, we consider a one-dimensional linear diffusion with exponential quadratic cost functional. Perhaps, the easiest example of a linear-quadratic (LQ) risk-sensitive control problem with mean-field coupling is 
\begin{eqnarray}
\label{LQ} 
\left\{
\begin{array}{lll} \inf_{u(\cdot)\in \cU}  E^u e^{\theta \left[ \frac{1}{2} \int_0^T u^2(t)dt +\frac{1}{2}x^2(T)+\mu  E^u [x(T)] \right]},
\\
\displaystyle{\mbox{ subject to }\ }\\
dx(t)=\left(ax(t)+bu(t)\right)dt+\sigma dW_t+\alpha d\widetilde W_t,\\
dY_t=\b x(t)dt+d\widetilde W_t
\\
x(0)=x_{0},\, Y_0=0,\\
 \end{array}
\right.
\end{eqnarray}
where, $a, b, \a,\b,\mu$ and $\sigma$ are real constants.

\ms\no In this section we will illustrate our approach by only considering the LQ risk-sensitive control under partial observation without the mean-field coupling i.e. $(\mu=0)$ so that our result can be compared with \cite{charala} where a similar example (in many dimensions) is studied using the Dynamic Programming Principle. The case $\mu\neq 0$ can treated in a similar fashion (cf. \cite{DT-1}).

\ms\no We consider the linear-quadratic risk-sensitive control problem:
\begin{eqnarray}
\label{LQ} 
\left\{
\begin{array}{lll} \inf_{u(\cdot)\in \cU}  E^u e^{\theta \left[ \frac{1}{2} \int_0^T u^2(t)dt +\frac{1}{2}x^2(T)\right]},
\\
\displaystyle{\mbox{ subject to }\ }\\
dx(t)=\left(ax(t)+bu(t)\right)dt+\sigma dW_t+\alpha d\widetilde W_t,\\
dY_t=\b x(t)dt+d\widetilde W_t
\\
x(0)=x_{0},\, Y_0=0,\\
 \end{array}
\right.
\end{eqnarray}
where, $a, b, \a,\b$ and $\sigma$ are real constants.

\ms\no An admissible process $(\bar\rho(\cdot), \bar x(\cdot), \bar u(\cdot))$ that satisfies the necessary optimality conditions of Theorem \ref{main result} is obtained by solving the following system of forward-backward SDEs (cf. (\ref{SDEu}) and (\ref{rs-pq})) (see Remark \ref{R-main}, above).  
\begin{equation}\label{SDEu-LQ}
\left\{\begin{array}{lll}
d\bar\rho(t)=\b\bar\rho(t)\bar x(t) dY_t,\\ 
d\bar x(t)=\left\{c\bar x(t)+b\bar u(t)\right\}dt+\sigma dW_{t} +\a dY_{t}, \\ 
dp(t)=-\left(\begin{array}{ccc} H^{\theta}_{\rho}(t) \\ H^{\theta}_x(t)\end{array}\right)dt+q(t)(-\theta \ell(t)dt+dB_t),\\
dv^{\theta}(t)=\theta v^{\theta}(t)\langle \ell(t),dB_t\rangle,\\
p(T)=-\left(\begin{array}{ccc} (\theta\bar\rho(T))^{-1}\\ \bar x(T) \end{array}\right),\\ v^{\theta}(T)=\psi^{\theta}_T,
\\ \bar\rho(0)=1,\,\, \bar x(0)=x_0.
\end{array}\right.
\end{equation}
where, 
$$ c:=a-\a\b,\,\, B_t:=\left(\begin{array}{lll} Y_t \\ W_t\end{array}\right),\,\, 
\ell:=\left(\begin{array}{ccc} \ell_1\\ \ell_2 \end{array}\right),\,\, p:=\left(\begin{array}{ccc} p_1\\ p_2 \end{array}\right),\,\, q:=\left(\begin{array}{ccc} q_{11} & q_{12}\\ q_{21} & q_{22} \end{array}\right),
$$

$$
\psi^{\theta}_T:=\bar \rho(T)e^{\theta \left[ \frac{1}{2} \int_0^T \bar u^2(t)dt +\frac{1}{2}\bar x^2(T)\right]},
$$
and the associated risk-sensitive Hamiltonian is 

\be\label{rs-hamiltonian-LQ}\begin{array}{lll}
H^{\theta}(t,\rho,x, u, p, q,\ell):=(cx+bu)p_2-\frac{1}{2}u^2+\rho\b x(q_{11}+\theta \ell_1 p_1)\\ \qquad\qquad\qquad\qquad\quad +\a(q_{21}+\theta\ell_2 p_1)+\sigma (q_{22}+\theta\ell_2 p_2).
\end{array}
\ee

\ms\no Below, we derive an explicit solution of the system (\ref{SDEu-LQ}) and characterize the optimal control of our problem. 
 
 \ms\no We have
 $$
 H^{\theta}_u=bp_2-u,\quad H^{\theta}_{\rho}=\b x(q_{11}+\theta \ell_1 p_1),\quad H^{\theta}_x=cp_2+\b \rho (q_{11}+\theta \ell_1 p_1).
 $$
 Therefore,  in view of Theorem \ref{main result}, if $\bar u$ is an optimal control of the system 
(\ref{LQ}), it is necessary that 
\be
E^{\theta}[bp_2(t)-\bar u(t)|\cF^Y_t]=0.
\ee
This yields 
\be\label{LQ-u}
\bar u(t)=bE^{\theta}[p_2(t)|\cF^Y_t].
\ee
The associated state dynamics $\bar x$ solves then the  SDE
\be\label{LQ-state}
d\bar x(t)=\left\{c\bar x(t)+b^2 E^{\theta}[p_2(t)|\cF^Y_t]\right\}dt+\sigma dW_{t} +\a dY_{t}
\ee
We try a solution of the form 
\be\label{LQ-p}
p_1(t):=-\lambda(t)/\bar\rho(t),\quad  p_2(t):=-\gamma(t) \bar x(t),
\ee
where, $\lambda(t)$ and $\gamma(t)$ are deterministic functions such that $\lambda(T)=1/\theta$ and $\gamma(T)=1$.

\ms\no Noting that $s(t):=\bar\rho^{-1}(t)$ satisfies the SDE
\be\label{LQ-s}
ds(t)=\b^2\bar x(t)^2s(t)dt-\b\bar x(t)s(t)dY_t,
\ee
we have
\be\label{LQ-p-x}
\left\{\begin{array}{lll}
dp_1(t)=-(\dot\lambda(t)+\b\bar x^2(t)\lambda(t))s(t)dt+\b\lambda(t)\bar x(t)s(t)dY_t,\\
dp_2(t)=-(\dot\gamma(t)\bar x(t)+c\gamma(t)+b\bar u(t))dt-\sigma\gamma(t)dW_t-\a\gamma(t)dY_t.
\end{array}\right.
\ee
Identifying the coefficients in (\ref{LQ-p-x}) with the corresponding ones in (\ref{SDEu-LQ}), we obtain
\be\label{LQ-q}
q_{11}(t)=\b\lambda(t)\bar x(t)s(t),\;\;\; q_{12}(t)=0,\;\;\; q_{21}(t)=-\a\gamma(t),\;\;\; q_{22}=-\sigma\gamma(t),
\ee
and
\be\label{LQ-lambda-gamma}\left\{\begin{array}{lll}
\dot\lambda(t)=0,\\
\left(\dot\gamma(t)+2c\gamma(t)-\b^2\lambda(t)\right)\bar x(t)+b\bar u(t)\gamma(t)+\theta\left(\b\lambda(t)+\a\gamma(t)\right)\ell_1(t)+\theta\sigma\gamma(t)\ell_2(t)=0,\\ \lambda(T)=1/\theta, \,\, \gamma(T)=1.
\end{array}\right.
\ee
Hence, 
$$
\lambda(t)=1/\theta,\quad 0\le t\le T,
$$ 
and  
\be
\left(\dot\gamma(t)+2c\gamma(t)-\b^2/\theta\right)\bar x(t)+b\bar u(t)\gamma(t)=-(\b+\theta\a\gamma(t))\ell_1(t)-\theta\sigma\gamma(t)\ell_2(t).
\ee
Therefore, in view of (\ref{LQ-u}) and (\ref{LQ-p}), we have
\be\label{LQ-gamma}
\left(\dot\gamma(t)+2c\gamma(t)-\b^2/\theta\right)\bar x(t)-b^2\bar \gamma^2(t)E^{\theta}[\bar x(t)|\cF^Y_t]=-(\b+\theta\a\gamma(t))\ell_1(t)-\theta\sigma\gamma(t)\ell_2(t).
\ee
Taking the conditional expectation, yields

\be\label{LQ-gamma}
\left(\dot\gamma(t)+2c\gamma(t)-\b^2/\theta-b^2\gamma^2(t)\right)E^{\theta}[\bar x(t)|\cF^Y_t]=-(\b+\theta\a\gamma(t))E^{\theta}[\ell_1(t)|\cF^Y_t]-\theta\sigma\gamma(t)E^{\theta}[\ell_1(t)|\cF^Y_t].
\ee
This equation is feasible only if we choose 
\be\label{LQ-ell-cond}
E^{\theta}[\ell_1(t)|\cF^Y_t]=\xi_1(t)E^{\theta}[\bar x(t)|\cF^Y_t],\quad E^{\theta}[\ell_2(t)|\cF^Y_t]=\xi_2(t)E^{\theta}[\bar x(t)|\cF^Y_t],
\ee
for some deterministic functions $\xi_1(t)$ and $\xi_2(t)$. This is possible if we choose e.g.
\be\label{LQ-ell}
\ell_1(t)=\xi_1(t)\bar x(t),\quad \ell_2(t)=\xi_2(t)\bar x(t).
\ee
In view of (\ref{SDEu-LQ}), the ansatz (\ref{LQ-ell}) makes the generic martingale $v^{\theta}$ satisfy the linear SDE
\be\label{LQ-v}
dv^{\theta}(t)=\theta v^{\theta}(t)\bar x(t)\left(\xi_1(t)dY_t+\xi_2(t)dW_t\right),
\ee 

\ms\no At this stage, the pair $(\xi_1,\xi_2)$  parametrizes a family of probability measures $P^{\theta}$ through $\ell$,  all equivalent with $P$, and which characterizes the optimal processes $(\rho(\cdot),\bar x(\cdot),\bar u(\cdot))$ through the SMP.

\ms \no Let us examine two typical cases (among many others).

\medskip {\bf Case 1}.  $\xi_1(t)=\xi_2(t)=1$. This choice yields the form $\ell(t)=(\bar x(t),\bar x(t))$, in the ansatz (\ref{LQ-ell}), which in turn gives a Riccati equation for $\gamma$:
\be\label{gamma-1}
\dot\gamma(t)+(2c+\theta(\a+\sigma))\gamma(t)-b^2\gamma^2(t)+\b-\b^2/\theta=0,\quad \gamma(T)=1,
\ee
whose solution is standard. 

\medskip {\bf Case 2}.  $\xi_1(t)=\xi_2(t)=\gamma(t)$. This choice yields the form $\ell(t)=\gamma(t)(\bar x(t),\bar x(t))$, in the ansatz (\ref{LQ-ell}),  which in turn gives a Riccati equation for $\gamma$:

\be\label{gamma-2}
\dot\gamma(t)+(2c+\b)\gamma(t)+(\theta(\a+\sigma)-b^2)\gamma^2(t)-\b^2/\theta=0,\quad \gamma(T)=1,
\ee
whose solution is also standard. 

\ms\no Given $\gamma$ which solves either (\ref{gamma-1}) (in which case $\ell_1(t)=\bar x(t)$) or (\ref{gamma-2}) (for which $\ell_1(t)=\gamma(t)\bar x(t)$), the corresponding optimal control is 
\be\label{LQ-u-1}
\bar u(t)=-b\gamma(t)E^{\theta}[\bar x(t)|\cF^Y_t],
\ee
where, in view of the filter equation displayed in Theorem 8.1 in \cite{liptser-shiryaev}, $\pi_t(\bar x):=E^{\theta}[\bar x(t)|\cF^Y_t]$ is solution of the SDE on $(\Omega,{\cal F},\dbF, \dbP^{\theta})$:
\be\label{LQ-hat-x}
\pi_t(\bar x)=x_0+\int_0^t (c-b^2\gamma(s))\pi_s(\bar x)ds+\a\int_0^t \left(1+\theta\left[\pi_s(\bar x\ell_1)-\pi_s(\bar x)\pi_s(\ell_1)\right]\right)d\bar Y^{\theta}_s,
\ee 
where, for $t\in[0,T]$, $\pi_t(\bar x\ell_1):=E^{\theta}[\bar x(t)\ell_1(t)|\cF^Y_t]$,  $\pi_t(\ell_1):=E^{\theta}[\ell_1(t)|\cF^Y_t]$ and  $\bar Y^{\theta}_t=Y_t-\theta\int_0^t \pi_s(\ell_1)ds$ is an $(\Omega,{\cal F},\dbF^Y, \dbP^{\theta})$-Brownian motion.


\bibliographystyle{plain}


\end{document}